\documentclass[11pt]{article}
\usepackage{amssymb}
\usepackage{amsmath}
\usepackage{amscd}
\usepackage{amstext}
\usepackage{amsfonts}
\usepackage{euscript}
\def\scr{\EuScript}
\usepackage{amssymb}

\newcommand{\De}{\underline D}

\newcommand{\od}{{\scr D}_{A/k}}
\newcommand{\odg}{{\scr D}_{\hat A /k}}
\newcommand{\N}{\mathbb {N}}

\newcommand{\pd}[1]{\frac{\partial}{\partial X_{#1}}}
\newcommand\Del{\underline{\Delta}}
\newcommand{\X}{{\bf X}}
\newcommand{\T}{{\bf T}}

\newcommand{\Sum}[1]{{\displaystyle \sum #1}}

\newcommand{\G}{\mathfrak}
\newcommand\GD{\underline{{\G D}}}

\newcounter{numero}[subsection]
\renewcommand{\thenumero}{(\thesection .\arabic{numero})}
\newenvironment{corolario}{\medskip
\refstepcounter{numero}\noindent {\sc  \thenumero\ Corollary.}\
\it}{\vspace{1ex}\par}

\newenvironment{teorema}{\medskip
\refstepcounter{numero}\noindent {\sc  \thenumero\ Theorem.}\
\it}{\vspace{1ex}\par}

\newenvironment{definicion}{\medskip
\refstepcounter{numero}\noindent {\sc  \thenumero\ Definition.}\
\it}{\vspace{1ex}\par}

\newenvironment{proposicion}{\medskip
\refstepcounter{numero}\noindent {\sc  \thenumero\ Proposition.}\
\it}{\vspace{1ex}\par}

\newenvironment{nota}{\medskip
\refstepcounter{numero}\noindent {\sc  \thenumero\ Remark.}\
}{\vspace{1ex}\par}

\newenvironment{demostracion}{
\noindent {\sc  Proof.}\ }{\hfill Q.E.D.\vspace{1ex}\par}

\newcommand{\numero}{\refstepcounter{numero}\noindent {\sc  \thenumero\ }}

\DeclareMathOperator{\HS}{HS} \DeclareMathOperator{\Der}{Der}
\DeclareMathOperator{\Ider}{IDer}

\DeclareMathOperator{\End}{End}

\DeclareMathOperator{\Frac}{Q}

\DeclareMathOperator{\rank}{rank}

\begin{document}

\title{ Hasse-Schmidt Derivations and Coefficient Fields in
Positive Characteristics}

\author{Magdalena Fern\'{a}ndez-Lebr\'{o}n and
Luis Narv\'{a}ez-Macarro\thanks{Both authors are partially supported
by DGI, BFM2001-3207.}\\ Department of Algebra, Faculty of
Mathematics
\\ University of Seville\\ E-mail: lebron@algebra.us.es,
narvaez@algebra.us.es }
\date{}
\maketitle

\thispagestyle{empty}

\begin{abstract}
We show how to express any Hasse-Schmidt derivation of an algebra
in terms of a finite number of them under natural hypothesis. As
an application, we obtain coefficient fields of the completion of
a regular local ring of positive characteristic in terms of
Hasse-Schmidt derivations.\medskip

\noindent Keywords: Regular local ring, Coefficient field,
Hasse--Schmidt derivation.\\ MSC 2000: 13H05, 13N15, 13N10, 13A35.
\end{abstract}

\section*{Introduction}

Let $k\xrightarrow{} A$ be a ring homomorphism. Hasse-Schmidt
derivations of $A$ over $k$ are generalizations of usual
derivations, but they do not carry an $A$-module structure.
Nevertheless, Hasse-Schmidt derivations have a non abelian group
structure lifting the addition of derivations.

In this paper we show how to express any Hasse-Schmidt derivation
in terms of a finite number of them under very reasonable
conditions. In proving our result, we find a natural way of
producing ``non-linear combinations" of Hasse-Schmidt derivations
which, to some extent, could play the role of the $A$-module
structure of derivations.

As an application, we express coefficient fields of the completion
of a regular local ring of positive characteristic in terms of
Hasse-Schmidt derivations, generalizing a similar result in
characteristic zero.
\medskip

Let us now comment on the content of this paper.
\medskip

In section 1 we recall the notions of Hasse-Schmidt derivation and
differential operator.
\medskip

Section 2 deals with the main result of this paper, namely that
any Hasse-Schmidt derivation $\GD$ can be expressed as
``non-linear combination" of a finite number of them  $\De^1,
\dots, \De^n$ whenever their degree 1 components
$D^1_1,\dots,D^n_1$ generate the module of usual derivations.
\medskip

In section 3 we apply our main result to generalize a well known
theorem of Nomura and to obtain coefficient fields of the
completion of noetherian local regular rings in the positive
characteristics case.
\medskip

Our results seem related to some results in \cite{traves-2000}. We
hope to return to this relationship in a future work.

We thank Herwig Hauser for pointing out a gap in the statement of
proposition \ref{main} in an earlier version of this work.

\section{Preliminaries and notations}

 All rings and algebras considered in this paper are assumed to
be commutative with unit element.
\medskip

\numero {\it Hasse-Schmidt derivations} (cf. \cite{has37} and
\cite{mat_86}, \S 27).
\medskip

Let $k\xrightarrow{f} A\xrightarrow{g} B$ be ring homomorphisms.
Let $t$ be an indeterminate over $B$, and set $B_m=B[t]/(t^{m+1})$
for $m\geq 0$ and $B_{\infty}=B[[t]]$. We can view $B_m$ as a
$k$--algebra in a natural way (for $m\leq \infty$).

A {\em Hasse--Schmidt derivation} (over $k$) of length $m\geq 1$
(resp. of length $\infty$) from $A$ to $B$, is a sequence
$\underline{D}=(D_0, D_1,\dots , D_m)$ (resp. $\underline{D}=(D_0,
D_1,\dots )$) of $k$--linear maps $D_i:A \longrightarrow B$,
satisfying the conditions: $$ D_0=g, \quad
D_i(xy)=\sum_{r+s=i}D_r(x)D_s(y) $$ for all $x,y \in A$ and all
$i>0$. In particular, the first component $D_1$ is a
$k$-derivation from $A$ to $B$. Moreover, $D_i$ vanishes on $f(k)$
for all $i>0$.

Any Hasse-Schmidt derivation $\underline{D}$ is determined by a
ring homomorphism $$E: x \in A \mapsto \sum_{i=0}^m D_i(x)t^i \in
B_m$$ with $E(x)\equiv x \mod t$.

When $A=B$ and $g=1_A$, we simply say that $\underline{D}$ is a
Hasse-Schmidt derivation of $A$ (over $k$). We write
$\HS_k(A,B;m)$ for the set of all Hasse-Schmidt derivations (over
$k$) of length $m$ from $A$ to $B$,
$\HS_k(A,B)=\HS_k(A,B;\infty)$, $\HS_k(A;m) = \HS_k(A,A;m)$ and
$\HS_k(A)=\HS_k(A,A;\infty )$.
\medskip

We say that a $k$-derivation $\delta:A\to B$ is {\em integrable}
\cite{mat-intder-I} if there is a Hasse-Schmidt derivation
$\underline{D}\in \HS_k(A,B)$ such that $D_1=\delta$. The set of
integrable $k$-derivations from $A$ to $B$, denoted by
$\Ider_k(A,B)$, is a submodule of the $k$-derivations $B$-module
$\Der_k(A,B)$.
\bigskip

\numero {\it Differential operators} (cf. \cite{ega_iv_4}, $\S
16$, $16.8$).
\medskip

Let $f:k\to A$ be a ring homorphism.
\medskip

For all $i\geq 0$, we inductively define the subsets
$\od^{(i)}\subseteq \End_k(A)$ in the follo\-wing way:
$$\od^{(0)}:=A\subseteq End_k(A),\quad\od^{(i+1)}:=\{ \varphi \in
End_k(A) \ | \ [\varphi,a]\in \od^{(i)},\ \ \forall a\in A \}.$$
The elements of $\od:={\displaystyle \bigcup_{i\geq 0} \od^{(i)}}$
(resp. of $\od^{(i)}$) are called {\em linear differential
operators} (resp. {\em linear differential operators of order
$\leq i$}) of $A/k$. The family
 $\{ \od^{(i)} \}_{i\geq 0}$ is an increa\-sing sequence of
$(A,A)$--bimodules of $\End_k(A)$ satisfying: $$\od^{(1)}=A\oplus
Der_k(A),\quad \od^{(i)} \circ \od^{(j)}\subset \od^{(i+j)},$$ and
$[P,Q]\in \od^{(i+j-1)}$ for all $P\in \od^{(i)}$, $Q\in
\od^{(j)}$. Hence, $\od$ is a filtered subring of $\End_k(A)$.
Moreover, linear differential operators of $A/k$ are
$I$--continuous for any $I$--adic topology. In particular, for any
linear differential operator $P$ of $A/k$, there is a unique
extension $\widehat P\in \odg$ to the completion  $\widehat A$ of
$A$ for some separated $I$--adic topology.

For each $\De\in \HS_k(A;m)$, one easily proves that $D_i \in
\od^{(i)}$ and then there is a unique extension
$\underline{\widehat{D}} \in \HS_k(\widehat{A};m)$.

In a similar way, if $S\subset A$ is a multiplicatively closed
subset, any Hasse-Schmidt derivation of $A/k$ extends uniquely to
a Hasse-Schmidt derivation of $S^{-1}A/k$.
\bigskip

\numero {\it Taylor expansions} (cf. \cite{mou-vi}).
\medskip

Let $n\geq 1$ be an integer. We write $\X =(X_1,\dots ,X_n)$, $\T
=(T_1,\dots, T_n)$, $\X +\T =(X_1+T_1,\dots,X_n+T_n)$ and, for
$\alpha,\beta\in \N^n$,  $\X^{\alpha} =X_1^{\alpha_1} \cdots
X_n^{\alpha_n}$, $|\alpha|=\alpha_1 +\cdots +\alpha_n $, $\alpha
!=\alpha_1 !\cdots \alpha_n !$ and $\binom{\beta}
{\alpha}=\binom{\beta_1}{\alpha_1} \cdots
\binom{\beta_n}{\alpha_n}$.

We consider the usual partial ordering in  $\N^n $: $\beta \geq
\alpha$ means $\beta_1 \geq \alpha_1, \dots ,\beta_n \geq
\alpha_n$. We write $\beta >\alpha $ if $\beta \geq \alpha$ and
$\beta \neq \alpha$.

Let $A$ be the formal power series ring $k[[\X ]]$ (or the
polynomial ring $A=k[\X ]$). For any $f(\X )=\Sum{_{\alpha \in
\N^n}} \lambda_{\alpha} \X^{\alpha} \in A$ we
 define $\Delta^{(\alpha)}(f(\X ))$ by: $\displaystyle
f(\X +\T )=\Sum{_{\alpha \in \N^n }} \Delta^{(\alpha )}(f(\X ))
\T^{\alpha}$. It is well known that (cf. \cite{ega_iv_4}, \S 16,
16.11): $\Delta^{(\alpha )}  \in \od^{(|\alpha|)}$,
$\Delta^{(\alpha )} (f\cdot g)=\Sum{_{\beta +\sigma =\alpha}}
\Delta^{(\beta)}(f) \Delta^{(\sigma)}(g)$, $\alpha
!\Delta^{(\alpha )}=(\frac{\partial}{\partial
X_1})^{\alpha_1}\cdots (\frac{\partial}{\partial X_n})^{\alpha_n}
$ and $\od^{(i)}={\displaystyle\bigoplus_{|\alpha|\leq i} A\cdot
\Delta^{(\alpha)}}= {\displaystyle\bigoplus_{|\alpha|\leq i}
\Delta^{(\alpha)}\cdot A}$.
\medskip

\numero \label{tay} If we denote $\Delta^{(0,\dots ,\stackrel{(j)}
i ,\dots ,0)}=\Delta_i^j$, then
$\underline{\Delta}^j=(1_A,\Delta_1^j, \Delta_2^j,\dots ) \in
\HS_k(A)$.
\medskip

Finally, let us recall the notion of quasi-coefficient field of a
local ring.
\medskip

\begin{definicion} {\em (cf. \cite{mat_80}, $38.F$)}
Let $(R,{\G m})$ be a local ring, $K=R/{\G m}$ and $k_0$ a
subfield of $R$. We say that $k_0$ is a {\em quasi--coefficient
field} of $R$, if the extension $K/k_0$ is formally \'{e}tale {\em
(cf. \cite{ega_iv_4}, 17.1.1 and \cite{mat_80}, 38.E)}.
\end{definicion}

In the case of characteristic $0$, a field extension is formally
\'{e}tale if and only if it is separably algebraic {\em (cf.
\cite{mat_80}, $38.E$)}. On the other hand, any extension of
perfect fields of characteristic $p>0$ is formally \'etale.
\medskip

The following proposition is well known.
\medskip

\begin{proposicion}\label{etale}
Let $k_0 \xrightarrow{} k\xrightarrow{f} A\xrightarrow{g} B$ be
ring homomorphisms and let's suppose that the extension $k_0\to k$
is formally \'{e}tale. Then, $\HS_{k_0}(A,B;m)=\HS_k(A,B;m)$ for
any integer $m\geq 1$ or $m=\infty$.
\end{proposicion}

\section{Generating Hasse-Schmidt derivations}
\label{sec:HS}

 Throughout this section,
let $k\xrightarrow{f} A$ be a ring homomorphism.
\medskip

We consider the following partial ordering in $\N^n $: $ \beta
\succeq \alpha$ means that $\beta \geq \alpha$ and if $\alpha_i
=0$ then $\beta_i=0$. \medskip

Let $N\geq 2$ be an integer and $\GD,\De^1, \dots,
\De^n\in\HS_k(A;N)$. For each $\mu\in\N^n$ we write
$D_{\mu}=D_{\mu_1}^1 \circ \cdots \circ D_{\mu_n}^n$. Let $C_{ld}$
be elements in $A$, $1\leq d\leq n$, $1\leq l \leq N-1$, such that
\begin{equation}\label{grande}
 {\G D}_i = \sum_{m=1}^i \left( \sum_{
                         \begin{array}{c}
                         \scriptstyle |\lambda |=i,\ |\mu |=m \\
                        \scriptstyle  \lambda\succeq \mu
                         \end{array}}
\prod_{d=1}^n \sum_{
      \begin{array}{c}
      \scriptstyle l^d_1+ \cdots +l^d_{\mu_d}=\lambda_d\\
    \scriptstyle   l_q^d\geq 1
      \end{array}}
 \prod_{q=1}^{\mu_d} C_{l_q^dd}
\right) D_{\mu}
\end{equation}
for all $i=1,\dots,N-1$, where we write
 $$\sum_{
      \begin{array}{c}
    \scriptstyle  l^d_1+ \cdots +l^d_{\mu_d}=\lambda_d\\
\scriptstyle  l_q^d\geq 1
      \end{array}}
 \prod_{q=1}^{\mu_d} C_{l_q^dd}=
   1   \quad \text{\ if\ }\quad \mu_d=\lambda_d=0.$$

\begin{proposicion} \label{main} Under the above hypothesis,
the $k$-linear application $$ \delta = {\G D}_N - \sum_{m=2}^N
\left( \sum_{
                         \begin{array}{c}
                       \scriptstyle   |\lambda |=N,|\mu |=m\\
                       \scriptstyle   \lambda\succeq \mu
                         \end{array}}
\prod_{d=1}^n \sum_{
      \begin{array}{c}
  \scriptstyle  l^d_1+ \cdots +l^d_{\mu_d}=\lambda_d\\
   \scriptstyle     l_q^d\geq 1
      \end{array}}
 \prod_{q=1}^{\mu_d} C_{l_q^dd}\right) D_{\mu} $$
is a $k$-derivation of $A$.
\end{proposicion}

\begin{demostracion}
\noindent Let us take $a,b\in R$. Since
 $$D_{\mu} (a b)=
\Sum{_{\rho +\sigma =\mu}} D_{\rho}(a) D_{\sigma}(b)$$ we obtain
\begin{multline*}
%linea 1
\delta(ab) = {\G D}_N(a) b+a{\G D}_N(b)+ \sum_{\begin{array}{c}
\scriptstyle v+w=N\\ \scriptstyle 1\leq v,w\leq N-1
\end{array}}
{\G D}_v(a) {\G D}_w(b)-
\\
%linea 2
 - \sum_{m=2}^N \left( \sum_{ \begin{array}{c} \scriptstyle  |\lambda |=N,|\mu |=m \\
     \scriptstyle   \lambda\succeq \mu \end{array}    }
    \prod_{d=1}^n
     \sum_{   \begin{array}{c} \scriptstyle l^d_1+ \cdots +l^d_{\mu_d}=\lambda_d\\
     \scriptstyle      l_q^d\geq 1\end{array}    }
       \prod_{q=1}^{\mu_d} C_{l_q^dd} \right)
\cdot \left( \sum_{  \rho +\sigma=\mu} D_{\rho}(a) D_{\sigma}(b)
\right) =
\end{multline*}
\begin{multline*}
%linea 1
={\G D}_N(a)  b+a {\G D}_N(b) + \sum_{\begin{array}{c}
\scriptstyle v+w=N\\ \scriptstyle
 1\leq v,w\leq N-1 \end{array}}
       {\G D}_v(a) {\G D}_w(b)-
       \\
%linea 2
 - \sum_{m=2}^N
\left(     \sum_{    \begin{array}{c}  \scriptstyle
                          |\lambda |=N,|\mu |=m \\
                  \scriptstyle        \lambda\succeq \mu
                      \end{array}  }
 \prod_{d=1}^n    \sum_{       \begin{array}{c}
                          \scriptstyle       l^d_1+ \cdots +l^d_{\mu_d}=\lambda_d\\
                          \scriptstyle       l_q^d\geq 1
                               \end{array}          }
\prod_{q=1}^{\mu_d} C_{l_q^dd} \right) \cdot \left( D_{\mu}(a) b+a
D_{\mu}(b) \right) -
\\
 -\sum_{m=2}^N \left( \sum_{\begin{array}{c}
       \scriptstyle      |\lambda |=N,|\mu |=m \\
     \scriptstyle   \lambda\succeq \mu
             \end{array}}
\prod_{d=1}^n \sum_{\begin{array}{c}
   \scriptstyle    l^d_1+ \cdots +l^d_{\mu_d}=\lambda_d\\
   \scriptstyle   l_q^d\geq 1
       \end{array}}
 \prod_{q=1}^{\mu_d} C_{l_q^dd} \right) \cdot
\left( \sum_{\begin{array}{c}
   \scriptstyle     \rho +\sigma=\mu\\
    \scriptstyle    \rho ,\sigma >0
             \end{array}}
D_{\rho}(a) D_{\sigma} (b)\right) =
\\
= \delta(a)b+a\delta(b)  + \sum_{\begin{array}{c}  \scriptstyle
v+w=N\\    \scriptstyle
 1\leq v,w\leq N-1 \end{array}}
       {\G D}_v(a) {\G D}_w(b)-
 \end{multline*}
\begin{multline*}
%linea 1
- \sum_{m=2}^N \left( \sum_{\begin{array}{c}
     \scriptstyle        |\lambda |=N,|\mu |=m \\
     \scriptstyle   \lambda\succeq \mu
             \end{array}}
\prod_{d=1}^n \sum_{\begin{array}{c}
   \scriptstyle    l^d_1+ \cdots +l^d_{\mu_d}=\lambda_d\\
   \scriptstyle   l_q^d\geq 1
       \end{array}}
 \prod_{q=1}^{\mu_d} C_{l_q^dd} \right) \cdot
\left( \sum_{\begin{array}{c}
    \scriptstyle    \rho +\sigma=\mu\\
     \scriptstyle   \rho ,\sigma >0
             \end{array}}
D_{\rho}(a) D_{\sigma} (b)\right).
\end{multline*}
We need to prove that
\begin{multline*}
\sum_{\begin{array}{c}  \scriptstyle   v+w=N\\ \scriptstyle 1\leq
v,w\leq N-1 \end{array}}
       {\G D}_v(a) {\G D}_w(b)=\\
=\sum_{m=2}^N \left( \sum_{\begin{array}{c}
     \scriptstyle        |\lambda |=N,|\mu |=m \\
      \scriptstyle  \lambda\succeq \mu
             \end{array}}
\prod_{d=1}^n \sum_{\begin{array}{c}
    \scriptstyle   l^d_1+ \cdots +l^d_{\mu_d}=\lambda_d\\
    \scriptstyle  l_q^d\geq 1
       \end{array}}
 \prod_{q=1}^{\mu_d} C_{l_q^dd} \right) \cdot
\left( \sum_{\begin{array}{c}
   \scriptstyle     \rho +\sigma=\mu\\
     \scriptstyle   \rho ,\sigma >0
             \end{array}}
D_{\rho}(a) D_{\sigma} (b)\right),
\end{multline*}
but
\begin{multline*}
%linea 1
\sum_{\begin{array}{c}
         {\scriptstyle v+w=N}\\
         {\scriptstyle 1\leq v,w\leq N-1}\\
        \end{array}}
{\G D}_v(a) {\G D}_w(b)=\\
 \end{multline*}
\begin{multline*}
%linea 2
=\sum_{\begin{array}{c}
         {\scriptstyle v+w=N}\\
         {\scriptstyle 1\leq v,w\leq N-1}\\
        \end{array}}
 \left[ \sum_{r=1}^v \left( \sum_{\begin{array}{c}
            {\scriptstyle |\tau |=v,|\rho |=r} \\
            {\scriptstyle \tau\succeq \rho}\\
            \end{array}}
    {\displaystyle \prod_{d=1}^n}
         \sum_{\begin{array}{c}
                 {\scriptstyle l^d_1+ \cdots +l^d_{\rho_d}=\tau_d}\\
                 {\scriptstyle l_q^d\geq 1}\\
                \end{array}}
  {\displaystyle  \prod_{q=1}^{\rho_d}} C_{l_q^dd} \right)
  D_{\rho} (a)\right] \cdot  \\
%linea 3
 \cdot \left[ \sum_{s=1}^w
         \left( \sum_{\begin{array}{c}
                  {\scriptstyle |\omega |=w,|\sigma |=s} \\
                  {\scriptstyle \omega \succeq \sigma }\\
                 \end{array}}
     {\displaystyle \prod_{d=1}^n}
          \sum_{\begin{array}{c}
                 {\scriptstyle l^d_1+ \cdots +l^d_{\sigma_d}=\omega_d}\\
                 {\scriptstyle l_{q'}^d\geq 1} \\
                 \end{array}}
  {\displaystyle  \prod_{q'=1}^{\sigma_d}} C_{l_{q'}^dd} \right)
   D_{\sigma} (b) \right]=   \\
%linea 3
= \sum_{\begin{array}{c}
         {\scriptstyle v+w=N}\\
         {\scriptstyle 1\leq v,w\leq N-1}\\
        \end{array}}
 \sum_{m=2}^N
 \sum_{\begin{array}{c}
     {\scriptstyle r+s=m}\\
     {\scriptstyle 1\leq r \leq v}\\
     {\scriptstyle 1\leq s \leq w}\\
                       \end{array}}
 \sum_{\begin{array}{c}
         {\scriptstyle |\tau |=v,\ |\omega |=w} \\
         {\scriptstyle |\rho |=r,\ |\sigma |=s} \\
         {\scriptstyle \tau \succeq \rho\ \omega
                        \succeq \sigma}\\
        \end{array}}\\
%linea 4
  \left({\displaystyle \prod_{d=1}^n}
 \sum_{\begin{array}{c}
         {\scriptstyle l^d_1+ \cdots +l^d_{\rho_d}=\tau_d}\\
         {\scriptstyle l_q^d\geq 1} \\
        \end{array}}
  {\displaystyle  \prod_{q=1}^{\rho_d}} C_{l_q^dd}
  D_{\rho} (a)
\right) \cdot \left({\displaystyle \prod_{d=1}^n}
 \sum_{\begin{array}{c}
         {\scriptstyle l^d_1+ \cdots +l^d_{\sigma_d}=\omega_d}\\
         {\scriptstyle l_{q'}^d\geq 1} \\
        \end{array}}
  {\displaystyle  \prod_{q'=1}^{\sigma_d}} C_{l_{q'}^dd}
  D_{\sigma}(b)
\right)=
\end{multline*}
\begin{multline*}
%linea 5
=\sum_{m=2}^N \sum_{\begin{array}{c}
              {\scriptstyle |\lambda |=N,|\mu |=m} \\
                             {\scriptstyle \lambda\succeq \mu}
                              \end{array}}
\sum_{\begin{array}{c} \scriptstyle \rho +\sigma =\mu\\
\scriptstyle \rho,\sigma >0\end{array}} \sum_{  \begin{array}{c}
 {\scriptstyle \tau+\omega=\lambda} \\
                             {\scriptstyle \tau\succeq \rho,\omega\succeq \sigma}\\
                        \end{array} }\\
%linea 6
  \left({\displaystyle \prod_{d=1}^n}
 \sum_{\begin{array}{c}
         {\scriptstyle l^d_1+ \cdots +l^d_{\rho_d}=\tau_d}\\
         {\scriptstyle l_q^d\geq 1} \\
        \end{array}}
  {\displaystyle  \prod_{q=1}^{\rho_d}} C_{l_q^dd}
  D_{\rho} (a)
\right) \cdot \left({\displaystyle \prod_{d=1}^n}
 \sum_{\begin{array}{c}
         {\scriptstyle l^d_1+ \cdots +l^d_{\sigma_d}=\omega_d}\\
         {\scriptstyle l_{q'}^d\geq 1} \\
        \end{array}}
  {\displaystyle  \prod_{q'=1}^{\sigma_d}} C_{l_{q'}^dd}
  D_{\sigma}(b)
\right)=\\
%linea 7
 =\sum_{m=2}^N \sum_{\begin{array}{c}
 {\scriptstyle |\lambda |=N,|\mu |=m} \\
                             {\scriptstyle \lambda\succeq \mu}\\
                        \end{array}}
\sum_{\begin{array}{c} \scriptstyle \rho +\sigma =\mu\\
\scriptstyle \rho,\sigma >0\end{array}} \sum_{  \begin{array}{c}
 {\scriptstyle \tau+\omega=\lambda} \\
                             {\scriptstyle \tau\succeq \rho,\omega\succeq \sigma}\\
                        \end{array} } {\displaystyle \prod_{d=1}^n}\\
%linea 8
 \sum_{\begin{array}{c}
        {\scriptstyle l^d_1+ \cdots +l^d_{\rho_d}=\tau_d}\\
        {\scriptstyle l^d_{\rho_d+1}+ \cdots +l^d_{\mu_d}=\omega_d}\\
        {\scriptstyle l_q^d\geq 1}
       \end{array}}
\left(
 {\displaystyle  \prod_{q=1}^{\rho_d}} C_{l_q^dd} \cdot D_{\rho}(a) \right)
 \left( {\displaystyle
\prod_{q'=\rho_d+1}^{\mu_d}} C_{l_{q'}^dd} \cdot D_{\sigma}
(b)\right) =\\
%linea 9
 =\sum_{m=2}^N \sum_{\begin{array}{c}
              {\scriptstyle |\lambda |=N,|\mu |=m} \\
                             {\scriptstyle \lambda\succeq \mu}
                              \end{array}}
\sum_{\begin{array}{c} \scriptstyle \rho +\sigma =\mu\\
\scriptstyle \rho,\sigma >0\end{array}} {\displaystyle
\prod_{d=1}^n} \sum_{\begin{array}{c}
        {\scriptstyle l^d_1+ \cdots +l^d_{\mu_d}=\lambda_d}\\
        {\scriptstyle l_q^d\geq 1}
       \end{array}} \left(
{\displaystyle  \prod_{q=1}^{\mu_d}} C_{l_q^dd}  D_{\rho}(a)
D_{\sigma}(b) \right)=
 \end{multline*}
\begin{multline*}
%linea 10
=\sum_{m=2}^N \left( \sum_{\begin{array}{c}
                        {\scriptstyle |\lambda |=N,|\mu |=m} \\
                             {\scriptstyle \lambda\succeq \mu}
                                               \end{array}}
{\displaystyle \prod_{d=1}^n} \sum_{\begin{array}{c}
        {\scriptstyle l^d_1+ \cdots +l^d_{\mu_d}=\lambda_d}\\
        {\scriptstyle l_q^d\geq 1}
       \end{array}}
{\displaystyle  \prod_{q=1}^{\mu_d}} C_{l_q^d d} \right) \cdot
\left( \sum_{\begin{array}{c} \scriptstyle \rho +\sigma =\mu\\
\scriptstyle \rho,\sigma >0\end{array}  } D_{\rho}(a)
D_{\sigma}(b) \right)
\end{multline*}
and the proposition is proved.
\end{demostracion}

\begin{teorema}\label{largo} Let $m\geq 1$ be an integer or
$m=\infty$. Let $\De^1, \dots, \De^n\in \HS_k(A;m)$ be
Hasse-Schmidt derivations of $A/k$ such that their components of
degree 1, $D_1^1,\dots,D_1^n$, form a system of generators of the
$A$-module $\Der_k(A)$. Then, for any Hasse-Schmidt derivation
$\GD\in\HS_k(A;m)$
 there exist $C_{ld}\in A$, $1\leq d\leq n$, $1\leq l < m+1$,
such that the equation (\ref{grande}) holds for all $i\geq  1$.
Moreover, if $\{ D_1^1,\dots,D_1^n \}$ is a $A$-basis of
$\Der_k(A)$, then the $\{ C_{ld}\}$ are unique.
\end{teorema}
\begin{demostracion}
We proceed by induction on $i$.
\medskip

\noindent For $i=1$, ${\G D}_1$ is a derivation and so there exist
$C_{11},\dots ,C_{1n}\in A$ such that $${\G D}_1=C_{11}D^1_1
+\cdots +C_{1n}D^n_1.$$

\noindent Let $N$ be an integer $\geq 2$ and suppose we have
elements $C_{ld}\in A$, $1\leq l\leq N-1$, $1\leq d\leq n$ such
that relation (\ref{grande}) is true for $1\leq i\leq N-1$.
\medskip

By proposition \ref{main}, the $k$-linear application $$ \delta =
{\G D}_N - \sum_{m=2}^N \left( \sum_{
                         \begin{array}{c}
                       \scriptstyle   |\lambda |=N,|\mu |=m\\
                       \scriptstyle   \lambda\succeq \mu
                         \end{array}}
\prod_{d=1}^n \sum_{
      \begin{array}{c}
  \scriptstyle  l^d_1+ \cdots +l^d_{\mu_d}=\lambda_d\\
   \scriptstyle     l_q^d\geq 1
      \end{array}}
 \prod_{q=1}^{\mu_d} C_{l_q^dd}\right) D_{\mu} $$
is a $k$-derivation of $A$. Then, there exist $C_{Nd}\in A$,
$1\leq d\leq n$ s.t. $$ \delta = C_{N1}D_1^1 + \cdots +
C_{Nn}D_1^n$$ and
\begin{multline*}
%linea 1
{\G D}_N =\delta + \sum_{m=2}^N \left( \sum_{
                         \begin{array}{c}
                       \scriptstyle   |\lambda |=N,|\mu |=m\\
                       \scriptstyle   \lambda\succeq \mu
                         \end{array}}
\prod_{d=1}^n \sum_{
      \begin{array}{c}
  \scriptstyle  l^d_1+ \cdots +l^d_{\mu_d}=\lambda_d\\
   \scriptstyle     l_q^d\geq 1
      \end{array}}
 \prod_{q=1}^{\mu_d} C_{l_q^dd}\right) D_{\mu}= \\
 %linea 2
= \sum_{m=1}^N \left( \sum_{
                         \begin{array}{c}
                       \scriptstyle   |\lambda |=N,|\mu |=m\\
                       \scriptstyle   \lambda\succeq \mu
                         \end{array}}
\prod_{d=1}^n \sum_{
      \begin{array}{c}
  \scriptstyle  l^d_1+ \cdots +l^d_{\mu_d}=\lambda_d\\
   \scriptstyle     l_q^d\geq 1
      \end{array}}
 \prod_{q=1}^{\mu_d} C_{l_q^dd}\right) D_{\mu}
 \end{multline*}
and equation (\ref{grande}) holds for $i=N$.
\medskip

Obviously the $C_{ld}$ are unique if $\{D_1^1,\dots,D_1^n\}$ is a
$A$-basis of $\Der_k(A)$.
\end{demostracion}

\begin{nota} The $C_{ld}$ in theorem \ref{largo} depend on the
order of the $\De^1, \dots, \De^n$.
\end{nota}

\section{Applications: Coefficient Fields in positive characteristics}

\numero \label{nume:1} Let $(R,{\G m},k)$ be a noetherian regular
local ring of dimension $n\geq 1$ containing a field, $X_1,\dots
,X_n\in {\G m}$ a regular system of parameters of $R$, $k_0\subset
R$ a quasi-coefficient field and $\widehat{R}$ the completion of
$R$. We can identify $\widehat{R} = k[[X_1,\dots,X_n]]$ by means
of a canonical $k_0$-isomorphism. Let $\pd{1} ,\dots , \pd{n}$ be
the usual basis of $\Der_k(\widehat{R})$ and $\Del^1,\dots
,\Del^n$ the Hasse--Schmidt derivations of $\widehat{R}$ over $k$
defined in \ref{tay}.
\medskip

Let us recall the following result of M. Nomura (\cite{mat-fps-I},
Th.~2.3, \cite{mat_86}, Th.~30.6)
\medskip

\begin{teorema}\label{99}
Under the hypothesis above, the following conditions are
equivalent:
\begin{enumerate}
\item[(1)] $\pd{i}(R)\subset R$ for all $i=1,\dots,n$.
\item[(2)] There exist $D_i\in \Der_{k_0}(R)$ and $a_i\in R$, $i=1,\dots,n$, such that $D_i(a_j)=\delta_{ij}$.
\item[(3)] There exist $D_i\in \Der_{k_0}(R)$ and $a_i\in R$, $i=1,\dots,n$ such
that $\det(D_i(a_j))\notin {\G m}$.
\item[(4)] $\Der_{k_0}(R)$ is a free $R$--module of rank $n$ (and
$\{ D_1,\dots ,D_n \}$ is a basis).
\item[(5)] $\rank \Der_{k_0}(R)=n$.
\end{enumerate}
\end{teorema}

The proof of the following corollaries are straightforward.
\medskip

\begin{corolario}\label{np} Under hypothesis and equivalent conditions of theorem \ref{99},
for any basis $D_1,\dots ,D_n \in \Der_{k_0}(R)$, their extensions
$\widehat{D}_1,\dots,\widehat{D}_n$ to $\widehat{R}$ form a basis
of $\Der_k(\widehat R)$. Moreover, the restrictions
$\pd{i}|_R:R\to R$, $i=1,\dots,n$, form a basis of
$\Der_{k_0}(R)$.
\end{corolario}

\begin{corolario}\label{car-0} Under hypothesis of corollary \ref{np},
let's suppose that $k_0$ is a field of characteristic $0$. Then,
the set $\{ a\in\widehat{R} \mid \widehat{D}_j(a)=0 \quad \forall
j=1,\dots,n\}$ is a coefficient field of $\widehat{R}$ (the only
one containing $k_0$).
\end{corolario}

The following theorem is an improvement of theorem \ref{99} and is
based on the results of section \ref{sec:HS}.
\medskip

\begin{teorema}\label{t99g}
Under the hypothesis of \ref{nume:1}, the following conditions are
equivalent:
\begin{enumerate}
\item[(1)] $\Delta_i^j (R)\subset R$, for all $j=1,\dots,n$, $i\geq 0$.
\item[(2)] There exist $ \De^i \in \HS_{k_0}(R) $
 and $a_i\in R$, $i=1,\dots,n$, such that $$D_i^j(a_l)=\left\{
\begin{array}{ll}
0 & i\geq 2, \ \forall j\\  \delta_{jl} & i=1,\ \forall j,l.
\end{array}
\right. $$
\item[(3)] There exist $D_i \in \Ider_{k_0}(R)$ and $a_i\in R$, $i=1,\dots,n$,
such that $\det(D_j(a_l))\notin {\G m}$.
\item[(4)] $\Der_{k_0}(R)$ is a free $R$--module of rank $n$ and
$\Ider_{k_0}(R)=\Der_{k_0}(R)$ (and $\{D_1,\dots,D_n\}$ is a
basis).
\item[(5)] $\rank \Ider_{k_0}(R)=n$.
\end{enumerate}
\end{teorema}

\begin{demostracion}

\noindent $(1)\Longrightarrow (2) \Longrightarrow (3)$, $(4)
\Longrightarrow (5) $ are straightforward.

\noindent $(3) \Longrightarrow (4) $ comes from theorem \ref{99}.
\medskip

\noindent $(5) \Longrightarrow (1)$: Let $\De^1, \dots ,\De^n \in
\HS_{k_0}(R)$ such that $D^1_1,\dots,D^n_1$ are linear independent
over $R$. Let us consider the extensions
$\widehat{\De}^1,\dots,\widehat{\De}^n\in \HS_{k}(\widehat{R})$,
whose degree 1 components $\widehat{D}^1_1,\dots,\widehat{D}^n_1$
are also linear independent over $\widehat{R}$.
\medskip

Following the lines of the proof of theorem \ref{largo}, we are
going to prove the following result: \smallskip

\noindent For any $j=1,\dots,n$, there exist $C_{ld}^j\in K =
\Frac (R)$, $1\leq d\leq n$, $1\leq l < +\infty$, such that
\begin{equation}\label{eq:grande-2}
 \Delta_i^j = \sum_{m=1}^i \left( \sum_{
                         \begin{array}{c}
                         \scriptstyle |\lambda |=i,\ |\mu |=m \\
                        \scriptstyle  \lambda\succeq \mu
                         \end{array}}
\prod_{d=1}^n \sum_{
      \begin{array}{c}
      \scriptstyle l^d_1+ \cdots +l^d_{\mu_d}=\lambda_d\\
    \scriptstyle   l_q^d\geq 1
      \end{array}}
 \prod_{q=1}^{\mu_d} C_{l_q^dd}^j
\right) \widehat{D}_{\mu}
\end{equation}
for all $i\geq 1$. \medskip

For\footnote{This is the same argument used in the proof of
theorem \ref{99}.}  $i=1$,  the $\Delta_1^j$ are $k$-derivations
of $\widehat{R}$ and then there exist $C_{1d}^j\in L=\Frac
(\widehat{R})$ such that $$\Delta_1^j = \sum_{d=1}^n C_{1d}^j
\widehat{D}^d_1.$$ Then, for any $m=1,\dots,n$ we have $$
\delta_{jm} = \Delta_1^j (X_m) = \sum_{d=1}^n C_{1d}^j
\widehat{D}^d_1(X_m)$$ and the matrix $(\widehat{D}^j_1(X_m))$
with entries in $R$ has a non-zero determinant. In particular
$C_{1d}^j\in K$.
\medskip

Let us suppose that for $N\geq 2$ there exist $C_{ld}^j\in K$,
$1\leq d,j\leq n$, $1\leq l \leq N-1$, such that
(\ref{eq:grande-2}) holds for $i=1,\dots, N-1$. We can consider
$\HS_k(\widehat{R})\subset \HS_k(L)$. By proposition \ref{main},
for any $j=1,\dots,n$ the $k$-linear application
\begin{equation} \label{eq:grande-3} \delta^j = \Delta_N^j - \sum_{m=2}^N \left( \sum_{
                         \begin{array}{c}
                       \scriptstyle   |\lambda |=N,|\mu |=m\\
                       \scriptstyle   \lambda\succeq \mu
                         \end{array}}
\prod_{d=1}^n \sum_{
      \begin{array}{c}
  \scriptstyle  l^d_1+ \cdots +l^d_{\mu_d}=\lambda_d\\
   \scriptstyle     l_q^d\geq 1
      \end{array}}
 \prod_{q=1}^{\mu_d} C_{l_q^dd}^j\right) \widehat{D}_{\mu}
 \end{equation}
is a $k$-derivation of $L$. Let $a\in R$ be a common denominator
for the $C_{ld}^j$, $1\leq d,j\leq n$, $1\leq l \leq N-1$. Then,
$$a^N \delta^j = a^N \Delta_N^j - \sum_{m=2}^N \left( \sum_{
                         \begin{array}{c}
                       \scriptstyle   |\lambda |=N,|\mu |=m\\
                       \scriptstyle   \lambda\succeq \mu
                         \end{array}}
\prod_{d=1}^n \sum_{
      \begin{array}{c}
  \scriptstyle  l^d_1+ \cdots +l^d_{\mu_d}=\lambda_d\\
   \scriptstyle     l_q^d\geq 1
      \end{array}}
 \prod_{q=1}^{\mu_d} (a^{l_q^d}C_{l_q^dd}^j)\right) \widehat{D}_{\mu}
$$ maps $\widehat{R}$ into $\widehat{R}$ and $a^N \delta^j\in
\Der_k(\widehat{R})$. There exist $\overline{C}_{Nd}^j\in L$,
$1\leq d,j\leq n$, such that $$ a^N \delta^j = \sum_{d=1}^n
\overline{C}_{Nd}^j \widehat{D}^d_1.$$Since the matrix
$(\widehat{D}^j_1(X_m))$ with entries in $R$ has a non-zero
determinant and $ (a^N \delta^j)(X_m)\in R$ (notice that
$\Delta_N^j(X_m)=0$), we deduce that $\overline{C}_{Nd}^j\in K$.
By setting $C_{Nd}^j = a^{-N}\overline{C}_{Nd}^j\in K$ we obtain
the expression (\ref{eq:grande-2}) for $i=N$.\medskip

From (\ref{eq:grande-2}) we deduce that $$\Delta_i^j (R)\subset
K\cap \widehat{R} = R,$$ for all $j=1,\dots,n$, $i\geq 0$, and (1)
is proved.
\end{demostracion}

\begin{nota} As noticed in \cite{mat_80}, page 289 for theorem \ref{99},
theorem \ref{t99g} also holds for $\HS_k(\widehat{R})\cap \HS(R)$
instead of $\HS_{k_0}(R)$, $\Der_k(\widehat{R})\cap \Der(R)$
instead of $\Der_{k_0}(R)$ and $\{\delta\in \Der(R)\ |\ \exists
\De\in \HS_k(\widehat{R})\cap \HS(R)\ \text{s.t.}\  D_1=\delta\}$
instead of $\Ider_{k_0}(R)$, and the mention to a
quasi-coefficient field can be avoided.
\end{nota}

\begin{nota} We do not know any example of a noetherian regular
local ring $(R,{\G m},k)$  containing  a quasi-coefficient field
$k_0$ (of positive characteristic) such that
$\Ider_{k_0}(R)\neq\Der_{k_0}(R)$.
\end{nota}

The following theorem generalizes corollary \ref{car-0} to the
case of  characteristic $p\geq 0$.
\medskip

\begin{teorema}\label{ct99g} Under the hypothesis of \ref{nume:1},
let $\De^1, \dots ,\De^n \in \HS_{k_0}(R)$ such that their degree
1 components $\{D^1_1,\dots,D^n_1\}$ form a basis of
$\Der_{k_0}(R)$. Let $\widehat{\De}^1,\dots,\widehat{\De}^n $ be
the extensions of $\De^1,\dots,\De^n$ to $\widehat{R}$. Then, the
set  $\{ a\in\widehat{R} \mid \widehat{D}_i^j(a)=0 \quad \forall
j=1,\dots,n ,\ i\geq 1 \}$ is a coefficient field of $\widehat{R}$
(the only one containing $k_0$).
\end{teorema}
\begin{demostracion} Since $\widehat{R}=k[[X_1,\dots
,X_n]]$, we have $k=\{ a\in\widehat{R} \mid \Delta_i^j(a)=0 \quad
j=1,\dots,n;\ i\geq 1\}$. By corollary \ref{np} we deduce that
$\{\widehat{D}^1_1,\dots,\widehat{D}_1^n\}$ is a
$\widehat{R}$--basis of $Der_k(\widehat{R})$, and from theorem
\ref{largo} we can express the $\Delta_i^j$ in terms of
$\widehat{D}_i^j$. In particular $$\{ a\in \widehat{R} \mid
\widehat{D}_i^j(a)=0 \quad \forall j=1,\dots,n ,\ i\geq 1 \}
\subset k.$$ \noindent The opposite inclusion comes from
proposition \ref{etale}.
\end{demostracion}


\begin{thebibliography}{00}

\bibitem{ega_iv_4}
A. Grothendieck and J. Dieudonn\'e:
\newblock ``\'El\'ements de G\'eom\'etrie Alg\'ebrique {IV}: \'Etude locale des
  sch\'emas et de morphismes de sch\'emas (Quatri\`eme Partie)", volume 32 of
  Inst. Hautes \'Etudes Sci. Publ. Math.,
\newblock Press Univ. de France, Paris, 1967.

\bibitem{has37}
 H. Hasse and F. K. Schmidt: Noch eine
Begr\"{u}ndung der Theorie der h\"{o}heren Differrentialquotienten
in einem algebraischen Funktionenkorper einer Unbestimmten.  {\em
J. Reine U. Angew. Math.} {\bf 177}(1937), 223-239.

\bibitem{mat-fps-I}
H. Matsumura:
\newblock Formal power series rings over polynomial rings {I}.
\newblock In ``Number theory, algebraic geometry and commutative algebra, in
  honor of Yasuo Akizuki", Kinokuniya, Tokyo,
  (1973), 511-520.

\bibitem{mat_80}
H. Matsumura:
\newblock ``Commutative algebra, second edition".
\newblock Benjamin/Cummings Publishing Co., Inc., Reading, Mass., 1980.

\bibitem{mat-intder-I}
H. Matsumura:
\newblock Integrable derivations.
\newblock  {\em Nagoya Math. J.} {\bf 87}(1982), 227--245.

\bibitem{mat_86}
H. Matsumura:
\newblock ``Commutative Ring Theory". Vol. 8 of Cambridge studies in
  advanced mathematics,
\newblock Cambridge Univ. Press, Cambidge, 1986.

\bibitem{mou-vi}
K. Mount and O. E. Villamayor:
\newblock Taylor series and higher derivations.
\newblock Universidad de Buenos Aires {\bf 18}, 1969 (reprinted 1979).

\bibitem{traves-2000}
W. N. Traves:
\newblock Tight closure and differential simplicity.
\newblock {\em J. Algebra} {\bf 228(2)}(2000), 457--476.


\end{thebibliography}
\end{document}